\documentclass[letterpaper, 10 pt, conference]{ieeeconf}
\IEEEoverridecommandlockouts                              
\usepackage{url}

\usepackage{amsfonts}
\usepackage{amssymb}
\usepackage{acronym}
\usepackage{color}
\usepackage[dvips]{graphicx}
\usepackage{graphicx}
\usepackage{epsfig}

\def\fskip#1{}

\newtheorem{theorem}{Theorem}

\newtheorem{assumption}{Assumption}

\newtheorem{corollary}{Corollary}

\newtheorem{lemma}{Lemma}

\renewcommand{\theenumi}{\arabic{enumi} )}

\def\E{\mathcal{E}}

\def\ol{\bar}

\def\b{\beta}

\def\a{\alpha}

\def\tl{\tilde}

\def\re{\mathbb{R}}

\def\dist{{\rm dist}}
\def\t{\tau}

\title{Distributed Subgradient Methods and Quantization Effects
\thanks{This research was partially supported by the
National Science Foundation under grant ECCS-0701623,
CAREER grants CMMI 07-42538  and DMI-0545910, and by DARPA ITMANET program.}}

\author{
\authorblockN{Angelia Nedi\'c}
\authorblockA{Department of Industrial and \\Enterprise
Systems Engineering\\University of
Illinois\\ Urbana-Champaign, IL 61801\\Email: angelia@uiuc.edu} \and
\authorblockN{Alex Olshevsky, Asuman Ozdaglar, and John N.\ Tsitsiklis}
\authorblockA{Department of Electrical Engineering and Computer
Science\\Laboratory for Information and Decision Systems\\
Massachusetts Institute of Technology\\
Cambridge, MA 02142\\
Email: \{alex\_o, asuman, jnt\}@mit.edu}}

\begin{document}
\maketitle

\thispagestyle{empty}
\pagestyle{empty}


\begin{abstract}
\noindent We consider a convex unconstrained optimization problem
that arises in a network of agents whose goal is to
cooperatively optimize the sum of the individual agent objective
functions through local computations and communications. For this
problem, we use averaging algorithms to develop distributed
subgradient methods that can operate over a time-varying topology.
Our focus is on the convergence rate of these methods and the
degradation in performance when only quantized information is
available. Based on our recent results on the convergence time of
distributed averaging algorithms, we derive improved upper bounds on
the convergence rate of the unquantized subgradient method. We then
propose  a distributed subgradient method under the additional
constraint that agents can only store and communicate quantized
information, and we provide bounds on its convergence rate that
highlight the dependence on the number of quantization levels.
\end{abstract}

\section{Introduction}

There has been much interest in developing distributed methods for
optimization in networked-systems consisting of multiple agents with
local information structures. Such problems arise in a variety of
environments including resource allocation among heterogeneous
agents in large-scale networks, and information  processing and
estimation in sensor networks. Optimization algorithms deployed in
such networks should be completely distributed, relying only on
local observations and information, and robust against changes in
network topology due to mobility or node failures.

Recent work \cite{distributed_subgrad} has proposed a subgradient
method for optimizing the sum of convex objective functions
corresponding to $n$ agents connected over a time-varying topology
(see also the short paper \cite{confsub}). The goal of the agents is
to cooperatively solve the unconstrained optimization
problem
\begin{eqnarray}
\begin{array}{cc}
\hbox{minimize } & \sum_{i=1}^n f_i(x)\cr
\hbox{subject to} & x\in\re^m, \label{problem}
\end{array}
\end{eqnarray}
where each $f_i:\re^m\to\re$ is a convex function, representing the
local objective function of agent $i$, and known only to this agent.  The
decision vector $x$ in problem (\ref{problem}) can be viewed as
either a resource vector whose components correspond to
resources allocated to each agent, or a global decision vector which
the agents are trying to optimize using local information. The
proposed method
builds on the work in
\cite{johnthes,
distasyn} (see also, \cite{distbook}). It
relies on every agent maintaining estimates of an
optimal solution to problem (\ref{problem}), and communicating these
estimates locally to its neighbors. The estimates are updated using
a combination of a subgradient iteration\footnote{For subgradient
methods see, for example, \cite{polyak, shor, polyakbook, hulem, nlp,
ourbook}.} and an averaging algorithm. The subgradient step
optimizes the local objective function while the averaging algorithm
is used to obtain information about the objective functions of the
other agents.

In this paper, we consider the distributed subgradient method
discussed in \cite{distributed_subgrad}, and provide improved
convergence rate results. In particular, we use our recent results
on the convergence time of averaging algorithms \cite{quant-con} and
establish new upper bounds on the difference between the objective
function value of the estimates of each agent and the optimal value
of problem (\ref{problem}). These bounds have a polynomial
dependence on the number of agents $n$ (in contrast with the error
bounds in \cite{distributed_subgrad, confsub}, which involve
exponential dependence on $n$). Furthermore, we study a variation of
the distributed subgradient method in which the agents have access
to quantized information, and provide bounds on the convergence time
that contain additional error terms due to quantization.

In addition to the papers cited above, our work is related to the
literature on reaching consensus on a particular scalar value or
on computing exact averages of the initial values of the agents,
a subject motivated by natural models of cooperative
behavior in
networked-systems (see, e.g., 
\cite{vicsek}, 
\cite{ali}, 
\cite{boyd}, 
\cite{reza},
\cite{spielman}, and 
\cite{alex_john_cdc,
alex_john}). Closely related is also the work in
\cite{quantized} and 
\cite{fagnani, carli_zampieri}, which study the effects of
quantization on the performance of averaging algorithms. Our work is
also related to the {\it utility maximization} framework for
resource allocation in networks
(see 
\cite{kelly}, 
\cite{low}, 
\cite{Srikant}). 
In contrast to this literature, however, we allow the local objective
functions to depend on the entire resource allocation vector.

The rest of this paper is organized as follows. In Section
\ref{sgd_method}, we describe the distributed subgradient method and
present an improved convergence rate estimate using our recently
established bounds on the convergence time of our averaging
algorithms \cite{quant-con}. In Section \ref{quantized_algo}, we
consider a version of the method under the additional constraint
that the agents can only exchange quantized information. We provide
convergence and rate of convergence results as a function of the
number of quantization levels. Section \ref{conclusions} contains
our concluding remarks.

\vskip .1pc

\noindent {\bf Notation and Basic Notions.}
We view all vectors as columns.
We use $e_i$ to denote the vector with $i$th entry equal to 1 and all
other entries equal to 0.
We use ${\bf 1}$ to denote a vector with all entries
equal to 1.
For a matrix $A$, we
use $a_{ij}$ or $[A]_{ij}$ to denote the matrix entry in the $i$th
row and $j$th column.
We write $[A]_i$ and $[A]^j$ to denote respectively the $i$th row
and the $j$th column of a matrix $A$.
A vector $a$ is said to be a {\it stochastic vector} when its
components $a_i$ are nonnegative and $\sum_i
a_i =1$. A square matrix $A$ is said to be {\it stochastic}
when each row of $A$ is a stochastic vector,
and it is said to be {\it doubly stochastic}
when both $A$ and its transpose $A'$ are stochastic matrices.

For a convex function $F:\re^m\to\re$, we use the
notion of a subgradient (see \cite{ourbook}):
a vector $s_F(\ol x)\in\re^m$ {\it is a
subgradient of a convex function $F$ at
$\ol x$} if
\[
F(\ol x) + s_F(\ol x)'(x-\ol x)\le F(x),\qquad \hbox{for all
}x. \label{sgdconvdef}
\]

We use the notation $f(x)=\sum_{j=1}^n f_i(x)$. We denote the
optimal value of problem (\ref{problem}) by $f^*$ and the
set of optimal solutions by $X^*$.

\section{Distributed Subgradient Method}\label{sgd_method}

We first introduce our distributed subgradient method for solving
problem (\ref{problem}) and discuss the assumptions imposed on the
information exchange among the agents. We consider a set
$V=\{1,\ldots,n\}$ of agents. Each agent starts with an initial
estimate $x_i(0)\in\re^m$ and updates its estimate at discrete times
$t_k, k=1,2,\ldots$. We denote by $x_{i}(k)$  the vector estimate
maintained by agent $i$ at time $t_k$. When updating, an agent $i$
combines its current estimate $x_i$ with the estimates $x_j$
received from its neighboring agents $j$. Specifically, agent $i$
updates its estimates by setting
\begin{equation}\label{gensyncrel}
x_i(k+1) = \sum_{j=1}^n a_{ij}(k)x_j(k) - \a d_i(k),
\end{equation}
where the scalars $a_{i1}(k),\ldots,a_{in}(k)$ are nonnegative
weights and the scalar $\a>0$ is a stepsize. The vector $d_i(k)$ is
a subgradient of the agent $i$ cost function $f_i(x)$ at $x=x_i(k)$.
 We use the notation $A(k)$ to denote the {\it weight matrix} $[a_{ij}(k)]_{i,j=1,\ldots,n}$.

The evolution of the estimates $x_i(k)$ generated by Eq.\
(\ref{gensyncrel}) can be equivalently represented using transition
matrices. In particular, we  define a {\it transition matrix}
$\Phi(k,s)$ for any $s$ and $k$ with $k\ge s$, as follows:
\[\Phi(k,s) = 
A(k)A(k-1)\cdots A(s+1)A(s) 
.\]
Using these transition matrices, we relate the estimate $x_i(k+1)$
to the estimates $x_1(s),\ldots,x_n(s)$ for any $s\le k$. In
particular, for the iterates generated by Eq.\ (\ref{gensyncrel}),
we have for any $i$, and any $s$ and $k$ with $k\ge s$,
\begin{eqnarray}
x_i(k+1)&=&\sum_{j=1}^n [\Phi(k,s)]_{ij} x_j(s)\cr & - & \a\,
\sum_{r=s}^{k-1} \sum_{j=1}^n [\Phi(k,r+1)]_{ij} d_j(r)\cr &-& \a\,
d_i(k) \label{trans_sys}
\end{eqnarray}
(for more details, see \cite{distributed_subgrad}). As seen from the
preceding relation, to study the asymptotic behavior of the
estimates $x^i(k)$, we need to understand the behavior of the
transition matrices $\Phi(k,s)$. We do this under some assumptions
on the agent interactions that translate  into some properties of
transition matrices.

Our first assumption imposes some conditions on the weights $a_{ij}(k)$ in Eq.\ (\ref{gensyncrel}).

\begin{assumption}
For all $k\ge 0$, the weight matrix $A(k)$
is doubly stochastic with positive diagonal.
Additionally, there is a scalar $\eta> 0$
such that if $a_{ij}(k)>0$, then $a_{ij}(k) \geq \eta$.
\label{weights}
\end{assumption}

The doubly stochasticity assumption on the weight matrix will
guarantee that the subgradient of the objective function $f_i$ of
every agent $i$ will receive the same weight in the long run. The
second part of the assumption states that each agent gives
significant weight to its own values and to the values of its
neighbors.

At each time $k$, the agents' connectivity can be represented by a directed
graph $G(k)=(V,\E(A(k)))$,
where $\E(A)$ is the set of directed edges $(j,i)$,
including self-edges $(i,i)$,
such that $a_{ij}>0$.
Our next assumption ensures that the agents are
connected frequently enough to persistently influence each other.

\begin{assumption} \label{connectivity}
There exists an integer $B \geq 1$
such that the directed graph
\[ \Big(V, \E(A(kB)) 
\bigcup \cdots \bigcup \E(A((k+1)B-1))
\Big)\]
is strongly connected for all $k\ge 0$.
\end{assumption}

\subsection{Preliminary Results}
Here, we provide some results that we use later in our convergence
analysis of method (\ref{gensyncrel}). These results hold
under Assumptions \ref{weights} and \ref{connectivity}.

Consider a related update rule of the form
\begin{equation}\label{scalarmix}
z(k+1) = A(k) z(k),
\end{equation}
where $z(0)\in\re^n$ is an initial vector\footnote{This update rule
captures the averaging part of Eq.\ (\ref{gensyncrel}), as it
operates on a particular component of the agent estimates, with the
vector $z(k)\in \re^n$ representing the estimates of the different agents
for that component.}. Define
\[V(k) = \sum_{j=1}^n (z_j(k) - \bar z(k) )^2
\qquad\hbox{for all }k\ge0,\] where $\bar z(k)$ is the average of
the entries of the vector $z(k)$. Under the doubly stochasticity of
$A(k)$, the initial average $\bar z(0)$ is preserved by the update
rule (\ref{scalarmix}), i.e., $\bar z(k)=\bar z(0)$ for all $k$.
Hence, the function $V(k)$ measures the ``disagreement'' in agent
values.

In the next lemma,  we give
a bound on the decrease of the agent disagreement
$V(kB)$, which is linear in $\eta$ and quadratic in $n^{-1}$.
This bound is an immediate consequence of Lemma 5 in \cite{quant-con}.

\vspace{0.1cm}
\begin{lemma} \label{lboundvar} Let Assumptions
\ref{weights} and \ref{connectivity} hold.
Then, $V(k)$ is nonincreasing in $k$. Furthermore, 
\[V((k+1)B)\le \left(1-\frac{\eta}{2 n^2}\right) V(kB)\quad\hbox{for all }
k\ge0.\]
\end{lemma}

\vspace {0.17cm}

Using Lemma \ref{lboundvar} we obtain the following result
for the transition
matrices $\Phi(k,s)$ of
Eq.\ (\ref{trans_sys}).

\begin{corollary} \label{convproduct}
Let Assumptions \ref{weights} and \ref{connectivity} hold.
Then, for all $i,j$ and all $k,s$ with $k\ge s$,
we have
\[
\left|[\Phi(k,s)]_{ij} - \frac{1}{n}\right| \leq \left(1 -
\frac{\eta}{4n^2}\right)^{\lceil{k-s+1\over B}\rceil-2}.\]
\end{corollary}

\vspace{0.15cm}
\begin{proof}
By Lemma \ref{lboundvar}, we have for all $k\ge s$,
\[V(kB)\le \left(1-{\eta\over 2n^2}\right)^{k-s}\, V(sB).\]
Let $k$ and $s$ be arbitrary with $k\ge s$, and let
$$\t B\le s <(\tau+1)B,\quad tB\le k <(t+1)B,$$
with $\tau\le t$.
Hence, by the nonincreasing property of $V(k)$,
we have
\begin{eqnarray*}
V(k)&\le& V(tB)\cr &\le&\left(1-{\eta\over 2n^2}\right)^{t-\tau-1}\,
V((\tau+1)B)\cr &\le& \left(1-{\eta\over 2n^2}\right)^{t-\tau-1}\,
V(s).
\end{eqnarray*}
Note that $k-s<(t-\t)B +B $ implying that ${k-s+1\over B}\le
t-\t+1$, where we used the fact that both sides of
the inequality are integers.
Therefore $\lceil{k-s+1\over B}\rceil-2\le t-\tau-1$, and
we have for all $k$ and $s$ with $k\ge s$,
\begin{equation}\label{vrel}
V(k)\le \left(1-{\eta\over 2n^2}\right)^{\lceil{k-s+1\over
B}\rceil-2} \,V(s).\end{equation} By Eq.\ (\ref{scalarmix}), we have
$z(k+1) = A(k) z(k)$, and therefore $z(k+1) = \Phi(k,s) z(s)$ for
all $k\ge s$. Letting $z(s)=e_i$ we obtain $z(k+1)=[\Phi(k,s)]^i$.
Using 
the inequalities (\ref{vrel})
and $V(e_i)\le 1$, we obtain
\[ V([\Phi(k,s)]^i) \leq
\left(1-{\eta\over 2n^2}\right)^{\lceil{k-s+1\over B}\rceil-2}.\]
The matrix $\Phi(k,s)$ is doubly stochastic, because it is the product
of doubly stochastic matrices. Thus, the average entry of
$[\Phi(k,s)]_i$ is $1/n$ implying that for all $i$ and $j$,
\begin{eqnarray*}
\left([\Phi(k,s)]_{ij} - \frac{1}{n}\right)^2 &\le& V([\Phi(k,s)]^i)
\cr &\le& \left(1-{\eta\over2n^2}\right)^{\lceil{k-s+1\over
B}\rceil-2}.
\end{eqnarray*}
From the preceding relation and
$\sqrt{1-\eta/(2n^2)}\le 1-\eta/(4n^2)$, we obtain
\[ \left|[\Phi(k,s)]_{ij} - \frac{1}{n}\right|
\le \left(1-{\eta\over 4n^2}\right)^{\lceil{k-s+1\over
B}\rceil-2}.\]
\end{proof}

\subsection{Convergence time}
We now study the convergence of the subgradient method
(\ref{gensyncrel}) and obtain a convergence time bound.  We assume
the uniform boundedness of the set of
subgradients of the cost functions $f_i$ at
all points\footnote{This assumption can be relaxed, see
\cite{distributed_subgrad}.}: for some scalar $L>0$, we have for all
$x\in \re^m$ and all $i$,
\begin{equation}\label{sgd_bounded}
\|g\|\le L\qquad\hbox{for all } g\in\partial f_i(x),
\end{equation}
where $\partial f_i(x)$ is the set of all subgradients
of $f_i$ at $x$.

We define the time-averaged
vectors $\hat x _i(k)$ of the iterates $x_i(k)$ generated by
Eq.\ (\ref{gensyncrel}), i.e.,
\begin{equation}\label{eq:avg} \hat{x}_i(k) = \frac{1}{k}
\sum_{h=0}^{k-1} x_i(h).
\end{equation}
The use of these vectors allows us to
bound the objective function improvement
at every iteration; see
\cite{distributed_subgrad, confsub}. Under the subgradient
boundedness assumption, we have the following result\footnote{The
assumption $\max_{1 \leq i \leq n} \|x_i(0)\| \leq \alpha L $ in
this theorem is not essential. We use this assumption mainly to
present a more compact expression for the bound on the convergence
time. A bound that explicitly depends on $\|x_i(0)\|$ can be
obtained by following a similar line of analysis}.

\begin{theorem}\label{convergresult}
Let Assumptions \ref{weights} and \ref{connectivity} hold, and
assume that the set $X^*$ of optimal solutions of
problem (\ref{problem}) is
nonempty. Let the sets of subgradients be bounded as in Eq.\
(\ref{sgd_bounded}). Also, let the initial vectors $x_i(0)$ in Eq.\
(\ref{gensyncrel}) be such that $\max_{1 \leq i \leq n} \|x_i(0)\|
\leq \alpha L $. Then, the averages $\hat x_i(k)$ of the
iterates obtained by the method (\ref{gensyncrel}) satisfy
\begin{eqnarray*}
f(\hat{x}_i(k))
& \le& f^*
+ \frac{n \mbox{ dist}^2(y(0),X^*)}{2\alpha k} \cr
&&+
\frac{\alpha L^2 C}{2} + 2 \alpha n L^2 C_1,
\end{eqnarray*}
where
\begin{equation}\label{newconstant}
C = 1 + 8 n C_1, \ C_1 = 1+{nB\over\beta(1-\beta)},
\ \beta=1-{\eta\over 4n^2},
\end{equation}
and $y(0)=(1/n)\sum_{i=1}^n x_i(0)$.
\end{theorem}

\vskip 0.2pc
\begin{proof}
The proof is identical to that of Proposition 3 in
\cite{distributed_subgrad} and relies on the use of our improved
convergence rate bound in Corollary \ref{convproduct}.
\end{proof}

The convergence rate result in the preceding theorem improves that of
Proposition 3 in \cite{distributed_subgrad},
where an analogous estimate is shown
with a worse value for the constant $\beta$.
In particular, there
the constant $\beta$ in  \cite{distributed_subgrad}
is given by $\b=1-\eta^{(n-1)B}$,
and $C_1$ increases exponentially with $n$.
As seen from Eq.\ (\ref{newconstant}), our
new constant $C_1$ increases only polynomially with $n$,
indicating a much more favorable scaling as the network size increases.

\section{Quantization effects}\label{quantized_algo}

We next study the effects of quantization on the convergence
properties of the subgradient method. In particular, we assume that
each agent receives and sends only quantized estimates, i.e.,
vectors whose entries are integer multiples of $1/Q$. At time $k$,
an agent receives quantized estimates $x_j^Q(k)$ from some of its
neighbors and updates according to the following rule:
\begin{equation}\label{quantizedsgd}
x_i^Q(k+1)= \left\lfloor\sum_{j=1}^n a_{ij}(k)x^Q_j(k)-\a \tl
d_i(k)\right\rfloor,
\end{equation}
where $\tl d_i(k)$ is a subgradient of $f_i$ at $x^Q_i(k)$, and
$\lfloor y\rfloor$ denotes the operation of (componentwise) rounding
the entries of a vector $y$ to the nearest multiple of ${1/ Q}$.
We also assume that the agents' initial estimates $x^Q_j(0)$ are
quantized.

To analyze the proposed method, we find it useful to rewrite Eq.\
(\ref{quantizedsgd}) as follows:
\begin{equation}\label{qsgderror}
x_i^Q(k+1)= \sum_{j=1}^n a_{ij}(k)x^Q_j(k)-\a \tl d_i(k) - e_i(k+1),
\end{equation}
where the error vector $e_i(k+1)$ is given by
\begin{equation}\label{agenterror}
e_i(k+1) = \sum_{j=1}^n a_{ij}(k)x^Q_j(k)-\a \tl d_i(k)
-x_i^Q(k+1).
\end{equation}
Thus, the method can be viewed as a subgradient method with external
(possibly persistent) noise, represented by $e_i(k+1)$.
Due to the rounding down to the nearest multiple of ${1/ Q}$,
the error vector $e_i(k+1)$ satisfies
\begin{equation}\label{errorentries}
0\le e_i(k+1)\le {1\over Q}\, {\rm\bf 1},\qquad \hbox{for all $i$ and
$k$},
\end{equation}
where the inequalities above hold componentwise.

Using the transition matrices $\Phi(k,s)$, we can rewrite the update equation
(\ref{qsgderror}) as
\begin{eqnarray}\label{xqlinear}
x_i^Q(k+1) & = & \sum_{j=1}^n [\Phi(k,0)]_{ij} x_j^Q(0) \cr && -\a
\sum_{s=1}^k\, \sum_{j=1}^n [\Phi(k,s)]_{ij}\tl d_j(s-1)\cr &&
-\sum_{s=1}^k\,\sum_{j=1}^n \,[\Phi(k,s)]_{ij}e_j(s)\cr &&-\a \tl
d_i(k) - e_i(k+1).
\end{eqnarray}
In addition, we
consider a related {\it stopped model}, where after some time $\tl
k$, {\it the agents cease computing subgradients $\tl d_j(k)$, and
also after time $\tl k+1$} stop quantizing  (so that they can
send and receive real numbers). Thus, in this stopped model, we have
$$\tl d_i(k)=0\hbox{ and }e_i(k+1)=0,
\quad\hbox{for all $i$ and $k\ge \tl k$}.$$

Let $\{\tl x_i(k)\}$, $i=1,\ldots,n$ be the sequences generated by
the stopped model,
associated with a particular time $\tl k$.
In view of the preceding relation, we have for
each $i$,
$$\tl x_i(k)=x^Q_i(k)\qquad
\hbox{for $k\le \tl k$},$$ and for $k\ge\tl k+1$,
\begin{eqnarray}\label{stoppedprocess}
\tl x_i(k) & = & \sum_{j=1}^n [\Phi(k,0)]_{ij} x_j^Q(0) \cr
&& -\a \sum_{s=1}^{\tl k}
\sum_{j=1}^n [\Phi(k,s)]_{ij} \tl d_j(s-1)\cr
&&-\sum_{s=1}^{\tl k} \sum_{j=1}^n [\Phi(k,s)]_{ij} e_j(s).
\end{eqnarray}
Using the result of Corollary \ref{convproduct}, we can show that the
stopped process converges as $k\to\infty$. In particular, we have
the following result.

\begin{lemma}\label{stoppedconv}
Let Assumptions \ref{weights} and \ref{connectivity} hold. Then, for
any $i$ and any $\tl k\ge0$, the sequence $\{\tl x_i(k)\}$ generated
by Eq.\ (\ref{stoppedprocess}) converges and the limit vector does
not depend on $i$, i.e.,
$$\lim_{k\to\infty} \tl x_i(k) = y(\tl k)\qquad
\hbox{for all $i$ and $\tl k$}.$$ Furthermore, for the limit
sequence $y(k)$, we have:
\begin{itemize}
\item [(a)] For all $k$,
$$y(k+1) = y(k) -{\a\over n}\, \sum_{j=1}^n \tl d_j(k)
-{1\over n}\,\sum_{j=1}^n e_j(k+1).\]
\item [(b)]
When the subgradient norms $\|\tl d_j(k)\|$ are uniformly bounded by
some scalar $L$ [cf.\ Eq.\ (\ref{sgd_bounded})] and the agents'
initial values are such that $\max_j\|x_j^Q(0)\|\le \a L$, then for
all $i$ and $k$, $\|x_i^Q(k)-y(k)\| \le$
\[2\,\left(\a L +{\sqrt{m}\over Q }\right)\,
\left(1 + {n B\over \b(1-\b)}\right),
\]
where $\beta = 1 - \frac{\eta}{4n^2}$ and $m$ is the dimension of
the vectors $x_i^Q$.
\end{itemize}
\end{lemma}
\vskip 0.3pc
\begin{proof}
By Corollary \ref{convproduct}, for any $s\ge0$, the entries
$[\Phi(k,s)]_{ij}$ converge to $1/n$, as $k\to\infty$. By letting
$k\to\infty$ in Eq.\ (\ref{stoppedprocess}), we see that the limit
$\lim_{k\to\infty} \tl x_i(k)$ exists and is independent of $i$.
Denote this limit by $y(\tl k)$, and note that it is given by
\begin{eqnarray}\label{ykseq}
 y(\tl k) &=&
{1\over n}\,\sum_{j=1}^n x_j^Q(0) -{\a\over n}\, \sum_{s=1}^{\tl k}
\sum_{j=1}^n \tl d_j(s-1) \cr
&&-{1\over n}\,\sum_{s=1}^{\tl k} \sum_{j=1}^n e_j(s).
\end{eqnarray}
From the preceding relation,
applied to different values of $\tl k$,
we see that
\[
y(k+1) = y(k) -{\a\over n}\, \sum_{j=1}^n \tl d_j(k) -{1\over
n}\,\sum_{j=1}^n e_j(k +1).
\] This establishes part (a).

Using the relations in Eqs.\ (\ref{xqlinear}) and
(\ref{ykseq}), and the subgradient boundedness,
we obtain for all $k$,
\begin{eqnarray*}
\|x_i^Q(k) \hspace{-8pt}&-& \hspace{-8pt}y(k)\|
\le   \sum_{j=1}^n
\left|[\Phi(k,0)]_{ij}-{1\over n}\right| \|x_j^Q(0)\| \cr
&&+ \a L
\sum_{s=1}^{k-1}\sum_{j=1}^n
\left|[\Phi(k,s)]_{ij}-{1\over n}\right|\cr
&&+\sum_{s=1}^{k-1}\sum_{j=1}^n
\left|[\Phi(k,s)]_{ij} -{1\over n}\right| \|e_j(s)\|\cr
&&+ 2\a L +\|e_i(k)\| +{1\over n}\sum_{j=1}^n \|e_j(k)\|.
\end{eqnarray*}
By using Corollary \ref{convproduct}, we
have for all $i$ and $j$,  and any $k\ge s$,
\begin{eqnarray*}
\|x_i^Q(k) \hspace{-8pt}& - &\hspace{-8pt}y(k)\|
\le  \sum_{j=1}^n \b^{\lceil{k+1\over
B}\rceil-2} \|x_j^Q(0)\| \cr
&+& \a L \sum_{s=1}^{k-1}\sum_{j=1}^n
\b^{\lceil{k-s+1\over B}\rceil-2}\cr
&+&
\sum_{s=1}^{k-1}\sum_{j=1}^n \b^{\lceil{k-s+1\over
B}\rceil-2}\|e_j(s)\| \cr
& +& 2\a L +\|e_i(k)\| +{1\over
n}\sum_{j=1}^n \|e_j(k)\|.
\end{eqnarray*}
Since $e_i(k)\le {\rm \bf 1}/Q$ [cf.\ Eq.\ (\ref{errorentries})], we
have
$$\|e_i(k)\|\le {\sqrt{m}\over Q}
\qquad\hbox{for all $i$ and $k$}.$$
From the preceding two relations, and the inequality
$\max_j \|x_j^Q(0)\|\le \a L$, we obtain for all $i$
and $k$,\\
\vskip 0.05pc\noindent $\|x_i^Q(k)-y(k)\| \le \a L n
\b^{\lceil{k+1\over B}\rceil-2}$
\begin{eqnarray*}
&&+\a L n\, \sum_{s=1}^{k-1} \b^{\lceil{k-s+1\over B}\rceil-2}
+{\sqrt{m}\over Q}\, n\, \sum_{s=1}^{k-1}\b^{\lceil{k-s+1\over
B}\rceil-2}\cr && +2\a L +2\,{\sqrt{m}\over Q}.
\end{eqnarray*}
By using
$\sum_{s=0}^{k-1} \b^{\lceil{k-s+1\over B}\rceil-2}\le
\frac{1}{\b}\,\sum_{r=0}^\infty \b^{\lceil{r+2\over B}\rceil-1}$, and
\[\sum_{r=0}^\infty \b^{\lceil{r+2\over B}\rceil-1}
=\sum_{r=0}^\infty \b^{\lceil{r+2\over B}\rceil-1}
\le B\sum_{t=0}^\infty \b^t
=\frac{B}{1-\beta},
\]
we finally obtain
\[
\|x_i^Q(k)-y(k)\| \le 2\left(\a L +{\sqrt{m}\over Q }\right)
\left(1 + {n B\over \b(1-\b)}\right).
\]
\end{proof}

According to part (a) of Lemma \ref{stoppedconv},
the vectors $y(k)$ can be viewed as the iterates
produced by the ``fictitious'' centralized algorithm:
\begin{equation}\label{ykiterate}
y(k+1) = y(k) -{\a\over n}\, \sum_{j=1}^n \tl d_j(k) -{1\over
n}\,\sum_{j=1}^n e_j(k+1),
\end{equation}
which is an approximate subgradient method with persistent
noise: The direction $\sum_{j=1}^n \tl d_j(k)$ is an approximate
subgradient of the objective function $f$
because each vector $\tl d_j(k)$ is a subgradient of
$f_i$ at $x_i^Q(k)$ instead of at $y(k)$. The error term
$(1/n)\sum_{j=1}^n e_j(k +1)$ can be viewed as the
noise experienced by the
whole system. The noise is persistent
since the magnitudes of the errors $e_j(k)$ are non-diminishing.

We now focus on establishing an error bound for the function values
at the points $y(k)$ of the stopped process of Eq.\
(\ref{ykiterate}), starting with $y(0)={1\over n} \sum_{j=1}^n
x_j^Q(0)$, and with the direction $\tl d_j(k)$ being a subgradient
of $f_j$ at $x_j^Q(k)$ for all $j$ and $k$. The process $y(k)$ is
similar to the stopped process analyzed in \cite{distributed_subgrad}, defined
using $x_j(k)$ in place of $x_j^Q(k)$. Thus, using the same analysis
as in \cite{distributed_subgrad} (see Lemma 5 therein),
we can show the following  basic result.

\begin{lemma} \label{basicrelfyk}
Let Assumptions \ref{weights} and \ref{connectivity} hold, and assume
that the set $X^*$ of optimal solutions of problem (\ref{problem}) is
nonempty. Let the sequence $\{y(k)\}$ be defined by
Eq.\ (\ref{ykiterate}), and the sequences $\{x_j^Q(k)\}$ for $j\in
\{1,\ldots,n\}$ be generated by the quantized subgradient method
(\ref{quantizedsgd}). Also, assume that the subgradients are
uniformly bounded as in Eq.\ (\ref{sgd_bounded}), and that
$\max_j\|x_j^Q(0)\|\le \a L$. Then,
the average vectors $\hat y(k)$ defined as in Eq.\ (\ref{eq:avg}),
satisfy for all $k\ge 1$,
\[f(\hat y(k))\le f^* +{n\, {\rm dist}^2(y(0),X^*)\over 2\a k} +
{\a L^2 \tl C\over 2},\]
where
$$\tl C=1+{8n \tl C_1\over \a L},$$
$$\tl C_1=\left(\a L +{\sqrt{m}\over Q }\right)
\left(1 + {n
B\over \b(1-\b)}\right),$$
$\b=1-{\eta\over 4n^2}$ and $y(0)={1\over n} \sum_{j=1}^n x_j^Q(0)$.
\end{lemma}

\vskip 0.3pc
\begin{proof}
Using the same line of analysis as in the proof of Lemma 5 in
\cite{distributed_subgrad}, we can show that for all $k$,\\
$\dist^2(y(k+1),X^*) \le  \dist^2(y(k),X^*)$
\begin{eqnarray*}
&+&{2\a\over n}
\sum_{j=1}^n\left(\|\tl d_j(k)\| + \|g_j(k)\|\right)
\|y(k)-x_j^Q(k)\|\cr
&-&{2\a\over n} \left[f(y(k))-f^*\right]
+{\a^2\over n^2}\sum_{j=1}^n\|\tl d_j(k)\|^2,
\end{eqnarray*}
where $g_j(k)$ is a subgradient of $f_j$ at $y(k)$.
By using the subgradient boundedness, we further obtain\\
$\dist^2(y(k+1),X^*) \le \dist^2(y(k),X^*)$
\begin{eqnarray*}
&& +{4\a L\over n} \sum_{j=1}^n \|y(k)-x_j^Q(k)\|\cr
&&-{2\a\over n} \left[f(y(k))-f^*\right] +{\a^2L^2\over n}.
\end{eqnarray*}
By using Lemma \ref{stoppedconv}(b), we have \\
\vskip 0.05pc\noindent
$\dist^2(y(k+1),X^*)
\le \dist^2(y(k),X^*)$
\[
+ 8\a L \tl C_1 -{2\a\over n}
\left[f(y(k))-f^*\right] +{\a^2L^2\over n},
\]
where
$\tl C_1=\left(\a L +{\sqrt{m}\over Q }\right)\, \left(1 + {n
B\over \b(1-\b)}\right)$.
Therefore, \\
\vskip 0.005pc\noindent
$f(y(k)) \le  f^* + {\a L^2\over 2} + 4nL\tl C_1$
\begin{eqnarray*}
&& + {n\over 2\a}\,
\left(\dist^2(y(k),X^*)-\dist^2(y(k+1),X^*)\right),
\end{eqnarray*}
and by regrouping the terms and introducing $\tl C=1+{8n \tl
C_1\over \a L}$, we have
for all $k$,\\
\vskip 0.005pc\noindent
$f(y(k)) \le f^* + {\a L^2 \tl C\over 2}$
\[ + {n\over 2\a}\,
\left(\dist^2(y(k),X^*)-\dist^2(y(k+1),X^*)\right).
\]
By adding these inequalities for different values of $k$, and by
using the convexity of $f$, we obtain the desired inequality.
\end{proof}

\vskip 0.3pc
Assuming that the agents can store real values
(infinitely many bits), we consider the time-average of the iterates
$x_i^Q(k)$, defined by
$${\hat x}^Q_i(k)={1\over k}\sum_{h=0}^{k-1}x_i^Q(h)
\qquad\hbox{for $k\ge 1$}.$$
Using Lemma \ref{basicrelfyk}, we have the following result.

\begin{theorem}\label{quantizconverg}
Under the same assumptions as in Lemma \ref{basicrelfyk},
the averages ${\hat
x}^Q_i(k)$ of the iterates obtained by the method
(\ref{quantizedsgd}) satisfy, for all $i$,
\[
f({\hat x}^Q_i(k)) \leq f^* + \frac{n\,{\rm dist}^2(y(0),X^*)}{2
\alpha k} + \frac{\alpha L^2 \tl C}{2} + 2 n L\tl C_1,
\]
where $\tl C$, $\tl C_1$, and $y(0)$ are as
in Lemma \ref{basicrelfyk}.
\end{theorem}

\vskip 0.3pc
\begin{proof}
By the convexity of the functions $f_j$, we have, for any $i$
and $k$,
$$f({\hat x}^Q_i(k))\le f(\hat y(k))
+\sum_{j=1}^n g_{ij}(k)'(\hat x^Q_i(k) -\hat y(k)),$$ where
$g_{ij}(k)$ is a subgradient of $f_j$ at $\hat x^Q_i(k)$. Then, by
using the boundedness of the subgradients and Lemma
\ref{stoppedconv}(b), we obtain for all $i$ and $k$,
\[
f({\hat x}^Q_i(k)) \le f(\hat y(k)) +2n L\tl C_1,\]
with $\tl
C_1=\left(\a L +{\sqrt{m}\over Q }\right)
\left(1 + {n B\over \b(1-\b)}\right)$. The result follows by using Lemma \ref{basicrelfyk}.
\end{proof}

When the quantization level $Q$
is increasingly finer (i.e., $Q\to
\infty$), the results of Theorems \ref{quantizconverg}
and~\ref{convergresult} coincide. More
specifically, when $Q\to \infty$, the constants $\tl C_1$ and $\tl C$
of Theorem \ref{quantizconverg} satisfy
$$\lim_{Q\to\infty}\tl C_1 = \a L
\left(1+{nB\over \b(1-\beta)}\right)= \a L C_1,$$
$$\lim_{Q\to\infty}\tl C
= 1 +\frac{8n}{\a L}\,
\lim_{Q\to\infty}\tl C_1= 1+ 8n C_1,$$
with
$C_1=1+{nB\over \b(1-\beta)}$. Thus, for the error term
of Theorem \ref{quantizconverg}, we have
\[\lim_{Q\to\infty}
\left(\frac{\alpha L^2 \tl C}{2} + 2 n L\tl C_1\right)
= \frac{\alpha L^2}{2} C
+ 2 n\a L^2 C_1
\]
where $C=1+ 8n C_1$ and
$C_1 = 1+{nB\over\b(1-\beta)}$.
Hence, in the limit as $Q\to\infty$,
the error terms in the estimate of Theorem
\ref{quantizconverg} reduce to the error terms
in the estimate of Theorem \ref{convergresult}.

\section{Conclusions}\label{conclusions}

We studied distributed subgradient methods for convex optimization
problems that arise in networks of agents connected through a
time-varying topology. We first considered an algorithm for the case
where agents can exchange and store continuous values, and proved a
bound on the convergence rate. We next studied the algorithm under
the additional constraint that agents can only send and receive
quantized values. We showed that our algorithm guarantees
convergence of the agent values to the optimal objective value
within some error. Our bound on the error highlights the dependence
on the number of quantization levels,
and the polynomial dependence on the number $n$ of agents.
Future work includes
investigation of the effects of other quantization schemes and of noise
in the agents' estimates on the performance of the algorithm.

\bibliographystyle{amsplain}
\bibliography{distributed_2}

\end{document}